\documentclass{amsart}
\usepackage{xcolor}
\usepackage{amsmath, amsfonts, amsthm}
\usepackage{mathrsfs}

\usepackage{hyperref}
\hypersetup{
    colorlinks=true,
    linkcolor=blue,
    filecolor=magenta,      
    urlcolor=cyan,
    citecolor=black
}
\newcommand{\Mod}[1]{\ (\mathrm{mod}\ #1)}

\newtheorem{theorem}{Theorem}[]

\newtheorem{proposition}[theorem]{Proposition}

\theoremstyle{remark}

\theoremstyle{definition}
\newtheorem{definition}[theorem]{Definition}

\newcommand{\Mrond}{\mathcal{M}}

\newcommand{\Orond}{\mathcal{O}}
\newcommand{\Prond}{\mathcal{P}}

\newcommand{\Crond}{\mathcal{C}}
\newcommand{\Srond}{\mathcal{S}}

\newcommand{\Rrond}{\mathcal{R}}

\newcommand{\Frond}{\mathcal{F}}

\newcommand{\Irond}{\mathcal{I}}

\newcommand{\PP}{\mathbb{P}}

\newcommand{\ZZ}{\mathbb{Z}}

\newcommand{\Nbf}{\mathbf{N}}

\newcommand{\rarr}{\rightarrow}

\newcommand{\sg}{\Srond_g}
\newcommand{\sgbar}{\overline{\Srond}_g}
\newcommand{\mg}{\Mrond_g}
\newcommand{\mgbar}{\overline{\Mrond}_g}

\newcommand{\mgnbar}{\overline{\Mrond}_{g, n}}
\newcommand{\sgn}{\Srond_{g, 2n}}
\newcommand{\sgnbar}{\overline{\Srond}_{g, 2n}}

\newcommand{\rgn}{\Rrond_{g, 2n}}

\newcommand{\Pic}{\text{Pic}}

\newcommand{\nsrm}{\mathrm{ns}}

\begin{document}
\author{Bogdan-Petru Carasca}
\address{Institut für Mathematik, Humboldt Universität zu Berlin Unter den Linden 6, 10099 Berlin, Germany}
\email{bogdanpetru.carasca@gmail.com}

\begin{abstract}
    The moduli space $\sgn$ parametrizes pointed curves with spin structure. We prove that $\Srond_{2, 4}$, $\Srond_{2, 6}$, $\Srond_{3, 2}$, $\Srond_{3, 4}$, $\Srond_{3, 6}$, $\Srond_{4, 2}$, $\Srond_{4, 4}$, $\Srond_{5, 2}$ and $\Srond_{5, 4}$ are uniruled.
\end{abstract}

\title[Moduli of spin curves] {Uniruledness of some moduli spaces \\of pointed spin curves}

\maketitle
\let\thefootnote\relax

\section{Introduction}

The moduli space $\sgn$ of $2n$-pointed spin curves of genus $g$ parametrizes tuples $(X, x_1, \dots, x_{2n}, \eta)$ where $(X, x_1, \dots, x_{2n})$ is a $2n$-pointed smooth curve of genus $g$ and $\eta$ is a line bundle on $X$ such that $\eta^{\otimes 2} = \omega_X(-x_1 - \dots - x_{2n})$. The natural forgetful map $\pi_{g , 2n} : \sgn \rarr \Mrond_{g, 2n}$ is finite of degree $2^{2g}$.
\vskip 0.5em

These moduli spaces are in close connection with the classical moduli spaces of spin curves. A spin curve of genus $g$ is a pair $[X, \eta]$ where $\eta$ is a line bundle on $X$ such that $\eta^{\otimes 2} = \omega_X$, and the moduli space of spin curves of genus $g$ is denoted by $\Srond_g$. We say that a spin curve $[X, \eta]$ is even (resp. odd) if $h^0(X, \eta)$ is even (resp. odd). In the simultaneous papers \cite{mumford1971theta} and \cite{atiyah1971riemann} it was shown that even and odd spin curves do not mix. This means that given a family of spin curves over a connected base, the parity is constant. As a consequence of this, $\sg$ splits into two connected components $\sg^+$ and $\sg^-$, the moduli spaces of even and odd spin curves respectively. We denote by $\sgbar^+$ and $\sgbar^-$ the Cornalba compactifications of these moduli spaces, cf. \cite{cornalba1989moduli}.

\vskip 0.5em
Given an irreducible $n$-nodal curve $X = \widetilde{X}/\{x_i \sim y_i\}_{1 \leq i \leq n}$, a spin structure on $X$ is equivalent to a spin structure on $[\widetilde{X}, x_1, y_1, \dots, x_n, y_n]$ plus some gluing data. Therefore $\sgn$ occurs naturally in the boundary of $\overline{\Srond}_{g + n}^+$ and $\overline{\Srond}_{g + n}^-$, up to a finite factor given by the gluing data, and the problem of describing their birational geometry is quite natural.

\vskip 0.5em
The problem of the Kodaira classification of $\mg$ has a rich history. It was known from \cite{severi1915sulla} that $\mg$ is unirational for $g \leq 10$, and great progress has been made in recent years. Currently, we know from \cite{sernesi1981unirazionalita}, \cite{chang1984unirationality} and \cite{verra2005unirationality} that $\mg$ is unirational for $g \leq 14$ and that it is rationally connected for $g = 15$. In the other direction, by \cite{harris1982kodaira} and \cite{eisenbud1987kodaira}, $\mgbar$ is of general type for $g \geq 24$. In the recent paper \cite{farkas2020kodaira}, the long-lasting open cases of $g = 22$ and $g = 23$ have been shown to be of general type. Naturally, the same problems have been posed for $\mgnbar$, and Logan in \cite{logan2003kodaira} shows that all but finitely many $\mgnbar$ are of general type for $g \geq 4$. There remain, nonetheless, some open cases, and are especially worth mentioning the cases of $16 \leq g \leq 21$.

\vskip 0.5em
On the other hand, for moduli spaces of spin curves $\sgbar^\pm$, the Kodaira classification has been completed by 
Farkas and Verra in the papers \cite{farkas2010birational}, \cite{farkas2012moduli} and \cite{farkas2014geometry}. For $\sgbar^+$, it has been shown that it is uniruled for $g \leq 7$, of Kodaira dimension $0$ for $g = 8$ and of general type for $g \geq 9$, whereas $\sgbar^-$ is uniruled for $g \leq 11$ and of general type for $g \geq 12$.

\vskip 0.5em
This paper is concerned with the uniruledness of some of the moduli spaces of pointed spin curves $\sgn$. More precisely, we show the following:
\begin{theorem}
    The moduli spaces $\Srond_{2, 4}$, $\Srond_{2, 6}$, $\Srond_{3, 2}$, $\Srond_{3, 4}$, $\Srond_{3, 6}$, $\Srond_{4, 2}$, $\Srond_{4, 4}$, $\Srond_{5, 2}$ and $\Srond_{5, 4}$ are uniruled.    
\end{theorem}

\vskip 0.5em
The method we employ to show the uniruledness of $\sgn$ for small $g$ and $n$ relies on the results of Lelli-Chiesa, Knutsen and Verra in \cite{lelli2024uni}, which revolve around Nikulin surfaces of non-standard type. In that paper, the moduli spaces $\rgn$ of double covers of smooth curves of genus $g$ ramified at $2n$ points are shown to be (uni)rational for $2 \leq g \leq 5$ and $n$ sufficiently small.

\vskip 0.5em
To give context to the method of proof, it should be pointed out the curves on Nikulin surfaces of non-standard type have rather particular properties. For example, it is proved in \cite{ma2012unirationality} that curves on Nikulin surfaces of non-standard type have Clifford index at most 2, far from the generic case in high genus. It is therefore rather surprising that one is able to use Nikulin surfaces to tackle rationality questions about moduli of curves. This approach is nonetheless successful in our situation as we focus on moduli spaces of pointed spin curves of small genus.

\subsection*{Acknowledgements}

The research was partly conducted in the framework of the DFG-funded
research training group RTG 2965: From Geometry to Numbers, Project
number 512730679.

\vskip 0.5em
I would like to thank Prof. Gavril Farkas for suggesting this problem to me and all the valuable discussions and advice that he gave me along the way.

\section{Preliminaries}

Let us now recall some definitions.
\begin{definition}
    A \textit{Nikulin surface of genus $h$} is the datum of $(S, M, H)$ where $(S, H)$ is a $K3$ surface of genus $h$ and $M \in \Pic(S)$ is a line bundle on $S$ such that $2M \sim N_1 + \dots + N_8$ for eight disjoint $(-2)$-curves on $S$. Write
    \[
    \Nbf := \bigoplus \ZZ[N_i] \subset \Pic(S),
    \]
    for the \textit{Nikulin lattice} and set $\Lambda_h := \Nbf \oplus \ZZ[H]$.
    A Nikulin surface is of \textit{standard type} if the embedding $\Lambda_h \subset \Pic(S)$ is primitive, i.e. $\Pic(S) / \Lambda_h$ is torsion-free. Otherwise it is called of \textit{non-standard type}.
\end{definition}

It is shown in \cite{van2007nikulin} that if $(S, M, H)$ is a Nikulin surface of non-standard type, then the genus $h$ is odd and $\Lambda_h \subset \Pic(S)$ has index two.

\vskip 0.5em
Based on this observation, we introduce the following definition. Let $R$, $R' \in \Pic(S)$ be such that, up to reordering,
\[
R = \frac{H - N_1 - N_2}{2}, \quad R' = \frac{H - N_3 - \dots - N_8}{2} \quad \text{if } h \equiv 3 \Mod{4}
\]
and
\[
R = \frac{H - N_1 - N_2 - N_3 - N_4}{2}, \quad R' = \frac{H - N_5 - N_6 - N_7 - N_8}{2} \quad \text{if } h \equiv 1 \Mod{4}.
\]
Set $g := g(R)$ and $g' := g(R')$ so that by adjunction $g = \frac{h+1}{4}$ and $g' = \frac{h - 3}{4}$ if $h \equiv 3 \Mod{4}$, and $g = g' = \frac{h - 1}{4}$ if $h \equiv 1 \Mod{4}$.

\vskip 0.5em
We will retain the notation in \cite{lelli2024uni}. To recall, let $\Frond_h^{\Nbf, \nsrm}$ be the moduli space of Nikulin surfaces of non-standard type of genus $h$. A slight modification of this moduli space, introduced to address the fact that $R$ and $R'$ are numerically identical in the case $h \equiv 1 \Mod{4}$, is the moduli space 
\begin{align*}
\widehat{\Frond}_h^{\Nbf, \nsrm} := \bigl\{&(S, M, H, R) \colon (S, M, H) \in \Frond_h^{\Nbf, \nsrm} \text{ and }R \in \Pic(S) \\
&\text{ such that } H - 2R \text{ is a sum of }(-2)\text{-curves}\bigr\}
\end{align*}
It is clear that $\widehat{\Frond}_h^{\Nbf, \nsrm}$ is an étale double cover of $\Frond_h^{\Nbf, \nsrm}$.

\vskip 0.5em
Define the following projective bundles over these varieties. Over $\Frond^{\Nbf, \nsrm}_h$ we let 
\[
\Prond_g := \bigl\{(S, M, H, X) \colon (S, M, H) \in \Frond_h^{\Nbf, \nsrm} \text{ and }X \in |R|\bigr\}
\]
and 
\[
\Prond'_{g'} := \bigl\{(S, M, H, X') \colon (S, M, H) \in \Frond_h^{\Nbf, \nsrm} \text{ and } X' \in |R'|\bigr\},
\]
and over $\widehat{\Frond}_h^{\Nbf, \nsrm}$ we let
\[
\widehat{\Prond}_g := \bigl\{(S, M, H, R, X) \colon (S, M, H, R) \in \widehat{\Frond}_h^{\Nbf, \nsrm} \text{ and } X \in |R|\bigr\}.
\]

Finally, let $r_{g, 2} \colon \Prond_g \rarr \Rrond_{g, 2}$, $r'_{g', 6} \colon \Prond'_{g'} \rarr \Rrond_{g', 6}$ and $r_{g, 4} \colon \widehat{\Prond}_g \rarr \Rrond_{g, 4}$ be the maps defined as follows. Let $(S, M, H, X) \in \Prond_g$, and set $\{x_i\} = X \cap N_i$. The associated point in $\Rrond_{g, 2}$ is
\[
r_{g, 2}(S, M, H, X) = [X, x_1 + x_2, \Orond_X(M^\vee)].
\]
Similarly, for $(S, M, H, X') \in \Prond'_{g'}$, we set
\[
r'_{g', 6}(S, M, H, X') = [X, x_1 + x_2 + x_3 + x_4 + x_5, + x_6, \Orond_X(M^\vee)]
\]
and for $(S, M, H, R, X) \in \widehat{\Prond}_g$ set
\[
r_{g, 4}(S, M, H, R, X) = [X, x_1 + x_2 + x_3 + x_4, \Orond_X(M^\vee)].
\]

\begin{theorem}[cf. \cite{lelli2024uni}]
    The maps $r_{2, 4}$, $r'_{2, 6}$, $r_{3, 2}$, $r_{3, 4}$, $r'_{3, 6}$, $r_{4, 2}$, $r_{4, 4}$, $r_{5, 2}$ and $r_{5, 4}$ are dominant.
\end{theorem}

\subsection{Irreducibility of \texorpdfstring{$\sgn$}{}}

This subsection is devoted to proving the following result.
\begin{proposition}
    The moduli spaces of pointed spin curves $\sgn$ are irreducible for $n > 0$.
\end{proposition}

Before starting the proof, let us recall some definitions and notation regarding $\sgbar^\pm$. A \textit{quasi-stable $2n$-pointed curve} is a pointed curve $[X, x_1, \dots ,x_{2n}]$ obtained by blowing up some of the nodes of a nodal pointed curve at most once. The rational components arising in this way are called \textit{exceptional components}. A \textit{stable spin curve} of genus $g$ is a triple $[X, \eta, \beta]$ where $X$ is a quasi-stable curve of genus $g$, $\eta$ is a line bundle on $X$ of degree $g - 1$ such that $\eta|_E = \Orond_E(1)$ for every exceptional component $E$ of $X$, and $\beta : \eta^{\otimes 2} \rarr \omega_X$ is homomorphism of line bundles on $X$ which is generically non-zero on the non-exceptional components.

\vskip 0.5em
Finally, recall the divisors $A_0^\pm \subset \sgbar^\pm$, whose generic point is an irreducible even (resp. odd) spin curve $[\widetilde{X}/{x_1 \sim x_2}, \eta]$. They are irreducible by \cite{cornalba1989moduli}.

\begin{proof}
    We proceed by induction on $n$. If $n = 1$, consider the map from $A_0^+ \cup A_0^-$ to $\Srond_{g, 2/\ZZ_2}$ given by
    \[
    [X, \eta] \mapsto [\widetilde{X}, x_1 + x_2, \nu^*\eta(-x_1 - x_2)],
    \]
    where $\nu : \widetilde{X} \rarr X$ is the normalization map, $p \in X$ is the node and $\nu^{-1}(p) = \{x_1, x_2\}$. This is clearly surjective, and by \cite{cornalba1989moduli} (6.1), both restrictions to $A_0^+$ and $A_0^-$ are surjective. As these boundary divisors are irreducible, it follows that $\Srond_{g, 2/\ZZ_2}$ is irreducible as well. From this we can easily deduce the irreducibility of $\Srond_{g, 2}$. Let $\gamma : [0, 1] \rarr \Srond_{g, 2/\ZZ_2}$ be a non-contractible loop based at $[X, x_1 + x_2, \eta]$, $\gamma(t) = [X(t), x_1(t) + x_2(t), \eta_t]$. This can be lifted to a path $\gamma'(t)$ from $[X, x_1, x_2, \eta]$ to $[X, x_2, x_1, \eta]$ by continuously lifting $x_1(t) + x_2(t) \in X(t)_2$ to $X(t) \times X(t)$. Therefore, the monodromy of $\Srond_{g, 2} \rarr \Srond_{g, 2/\ZZ_2}$ is transitive, implying that $\Srond_{g, 2}$ is irreducible. 

    \vskip 0.5em
    As in \cite{cornalba1989moduli}, we can construct a compactification $\sgnbar$ of $\sgn$ by means of stable pointed spin curves. Consider the inclusion
    \begin{align*}
        \Crond^1{\Srond}_{g, 2n-2} :&= \overline{\Mrond}_{g, 1} \times_{\mgbar} \overline{\Srond}_{g, 2n-2} \subset \sgnbar, \\
        [X, x_1, \dots, x_{2n - 2}, t_1, \eta_X] &\mapsto [X \cup_{t_1 \sim q_1} E \cup_{t_2 \sim q_2} \PP^1, x_1, \dots, x_{2n}, \eta]
    \end{align*}
    where $\{[\PP^1, x_{2n - 1}, x_{2n}, t_2]\} = \Mrond_{0, 3}$, $E$ is the exceptional component, and $\eta$ is given by $\eta|_X = \eta_X$, $\eta|_E = \Orond_E(1)$ and $\eta|_{\PP^1} = \Orond(-2)$. Suppose that $\sgnbar$ is not irreducible. On one hand, for a generic point $[X, x_1, \dots, x_{2n - 2}, t_1, \eta_X] \in \Crond^1\sgnbar$ the corresponding point in $\sgnbar$ is smooth, on the other hand, $\Crond^1\overline{\Srond}_{g, 2n -2}$ is irreducible by induction, so it is contained in the intersection of all the irreducible components of $\sgnbar$. This implies that every point in the image is singular, a contradiction.
    \end{proof}

\section{Proof of the main theorem}
This section contains the proof of the main theorem. 
\begin{proof}[Proof of Theorem 1]
We proceed case by case. Throughout $[X, x_1, \dots, x_{2n}, \eta] \in \sgn$ is a generic spin curve and $\Orond_X(Z)$ is an odd theta characteristic on $X$ with corresponding effective divisor $Z$. Before going into the individual cases themselves, we make a couple of useful remarks. Let $P = \{C_t\}_{t \in \PP^1}$ be a pencil of curves on a K3 surface $S$. We note that, since K3 surfaces are not ruled, the pencil is not isotrivial. Assume now $Z = p_1 + \dots + p_{g -1}$ is in the base locus of $P$ such that $\Orond_{C_0}(Z)$ is a theta characteristic on $C_0$, i.e. $\Orond_{C_0}(2Z) \simeq \omega_{C_0}$. One also has by adjunction that $\Orond_S(C) \otimes \Orond_C \simeq \omega_C$, and it follows that $\Orond_C(Z)$ is a theta characteristic for every $C$ in the pencil. With these two remarks in mind, we proceed with the proof.

\vskip 0.5em
\underline{$g = 2$ and $2n = 4$}: Since the map $r_{2, 4} : \widehat{\Prond}_2 \rarr \Rrond_{2, 4}$ is dominant, there exists a non-standard Nikulin surface $(S, M, H, R) \in \widehat{\Frond}^{\Nbf, \nsrm}_9$ such that $X \in |R|$, $\{x_i\} = X \cap N_i$ and $\eta \otimes \Orond_X(-Z) = \Orond_X(M^\vee)$. Consider the linear system $|R \otimes \Irond_Z|$, and as $h^0(S, R) = g + 1$ and $\deg Z = g - 1$, it follows that there exists a pencil $P \subset |R \otimes \Irond_Z|$ containing $X$. The map
\[
P \dashrightarrow \Srond_{2, 4},\quad X' \mapsto [X', x_1', x_2', x_3', x_4', M^\vee \otimes \Orond_{X'}(Z)]
\]
where $\{x_i'\} = X' \cap N_i$ exhibits a rational curve through $[X, x_1, x_2, \eta]$.

\vskip 0.5em
\underline{$g = 2$ and $2n = 6$}: Since the map $r'_{2, 6} : \Prond'_2 \rarr \Rrond_{2, 6}$ is dominant, there exists a non-standard Nikulin surface $(S, M, H) \in \Frond^{\Nbf, \nsrm}_{11}$ such that $X \in |R'|$, $\{x_i\} = X \cap N_i$ and $\eta \otimes \Orond_X(-Z) = \Orond_X(M^\vee)$. Consider the linear system $|R' \otimes \Irond_Z|$, and as before there exists a pencil $P \subset |R' \otimes \Irond_Z|$ containing $X$. The map
\[
P \dashrightarrow \Srond_{2, 6},\quad X' \mapsto [X', x_1', x_2', x_3', x_4', x_5', x_6', M^\vee \otimes \Orond_{X'}(Z)]
\]
where $\{x_i'\} = X' \cap N_i$ gives a rational curve through $[X, x_1, x_2, x_3, x_4, x_5, x_6, \eta]$.

\vskip 0.5em
\underline{$g = 3$ and $2n = 2$}: Since the map $r_{3, 2} : \Prond_3 \rarr \Rrond_{3, 2}$ is dominant, there exists a non-standard Nikulin surface $(S, M, H) \in \Frond^{\Nbf, \nsrm}_{11}$ such that $X \in |R|$, $\{x_i\} = X \cap N_i$ and $\eta \otimes \Orond_X(-Z) = \Orond_X(M^\vee)$. Considering the linear system $|R \otimes \Irond_Z|$, there exists a pencil $P \subset |R \otimes \Irond_Z|$ which contains $X$. The map
\[
P \dashrightarrow \Srond_{2, 4},\quad X' \mapsto [X', x_1', x_2', M^\vee \otimes \Orond_{X'}(Z)]
\]
where $\{x_i'\} = X' \cap N_i$ exhibits a rational curve through $[X, x_1, x_2, \eta]$.

\vskip 0.5em
\underline{$g = 3$ and $2n = 4$}: Since the map $r_{3, 4} : \widehat{\Prond}_3 \rarr \Rrond_{3, 4}$ is dominant, there exists a non-standard Nikulin surface $(S, M, H, R) \in \widehat{\Frond}^{\Nbf, \nsrm}_{13}$ such that $X \in |R|$, $\{x_i\} = X \cap N_i$ and $\eta \otimes \Orond_X(-Z) = \Orond_X(M^\vee)$. As before, there exists a pencil $P \subset |R \otimes \Irond_Z|$ containing $X$. The map
\[
P \dashrightarrow \Srond_{3, 4},\quad X' \mapsto [X', x_1', x_2', x_3', x_4', M^\vee \otimes \Orond_{X'}(Z)]
\]
where $\{x_i'\} = X' \cap N_i$ exhibits a rational curve through $[X, x_1, x_2, x_3, x_4, \eta]$.

\vskip 0.5em
\underline{$g = 3$ and $2n = 6$}: Since the map $r'_{3, 6} : \Prond_3' \rarr \Rrond_{3, 6}$ is dominant, there exists a non-standard Nikulin surface $(S, M, H) \in \Frond^{\Nbf, \nsrm}_{15}$ such that $X \in |R'|$, $\{x_i\} = X \cap N_i$ and $\eta \otimes \Orond_X(-Z) = \Orond_X(M^\vee)$. Considering the linear system $|R' \otimes \Irond_Z|$, there exists a pencil $P \subset |R' \otimes \Irond_Z|$ containing $X$. The map
\[
P \dashrightarrow \Srond_{3, 6},\quad X' \mapsto [X', x_1', x_2', x_3', x_4', x_5', x_6', M^\vee \otimes \Orond_{X'}(Z)]
\]
where $\{x_i'\} = X' \cap N_i$ exhibits a rational curve through $[X, x_1, x_2, x_3, x_4, x_5, x_6, \eta]$.

\vskip 0.5em
\underline{$g = 4$ and $2n = 2$}: Since the map $r_{4, 2} : \Prond_4 \rarr \Rrond_{4, 2}$ is dominant, there exists a non-standard Nikulin surface $(S, M, H) \in \Frond^{\Nbf, \nsrm}_{15}$ such that $X \in |R|$, $\{x_i\} = X \cap N_i$ and $\eta \otimes \Orond_X(-Z) = \Orond_X(M^\vee)$. Considering again the linear system $|R \otimes \Irond_Z|$, there exists a pencil $P \subset |R \otimes \Irond_Z|$ containing $X$. The map
\[
P \dashrightarrow \Srond_{4, 2},\quad X' \mapsto [X', x_1', x_2', M^\vee \otimes \Orond_{X'}(Z)]
\]
where $\{x_i'\} = X' \cap N_i$ exhibits a rational curve through $[X, x_1, x_2, \eta]$.

\vskip 0.5em
\underline{$g = 4$ and $2n = 4$}: Since the map $r_{4, 4} : \widehat{\Prond}_4 \rarr \Rrond_{4, 4}$ is dominant, there exists a non-standard Nikulin surface $(S, M, H, R) \in \widehat{\Frond}^{\Nbf, \nsrm}_{17}$ such that $X \in |R|$, $\{x_i\} = X \cap N_i$ and $\eta \otimes \Orond_X(-Z) = \Orond_X(M^\vee)$. Looking at the linear system $|R \otimes \Irond_Z|$, there exists a pencil $P \subset |R \otimes \Irond_Z|$ containing $X$. The map
\[
P \dashrightarrow \Srond_{4, 4},\quad X' \mapsto [X', x_1', x_2', x_3', x_4', M^\vee \otimes \Orond_{X'}(Z)]
\]
where $\{x_i'\} = X' \cap N_i$ exhibits a rational curve through $[X, x_1, x_2, x_3, x_4, \eta]$.

\vskip 0.5em
\underline{$g = 5$ and $2n = 2$}: Since the map $r_{5, 2} : \Prond_5 \rarr \Rrond_{5, 2}$ is dominant, there exists a non-standard Nikulin surface $(S, M, H) \in \Frond^{\Nbf, \nsrm}_{19}$ such that $X \in |R|$, $\{x_i\} = X \cap N_i$ and $\eta \otimes \Orond_X(-Z) = \Orond_X(M^\vee)$. Considering the linear system $|R \otimes \Irond_Z|$, there exists a pencil $P \subset |R \otimes \Irond_Z|$ containing $X$. The map
\[
P \dashrightarrow \Srond_{2, 4},\quad X' \mapsto [X', x_1', x_2', M^\vee \otimes \Orond_{X'}(Z)]
\]
where $\{x_i'\} = X' \cap N_i$ exhibits a rational curve through $[X, x_1, x_2, \eta]$.

\vskip 0.5em
\underline{$g = 5$ and $2n = 4$}: Since the map $r_{5, 4} : \widehat{\Prond}_5 \rarr \Rrond_{5, 4}$ is dominant, there exists a non-standard Nikulin surface $(S, M, H, R) \in \widehat{\Frond}^{\Nbf, \nsrm}_{21}$ such that $X \in |R|$, $\{x_i\} = X \cap N_i$ and $\eta \otimes \Orond_X(-Z) = \Orond_X(M^\vee)$. Considering the linear system $|R \otimes \Irond_Z|$, there exists a pencil $P \subset |R \otimes \Irond_Z|$ containing $X$. The map
\[
P \dashrightarrow \Srond_{5, 4},\quad X' \mapsto [X', x_1', x_2', x_3', x_4', M^\vee \otimes \Orond_{X'}(Z)]
\]
where $\{x_i'\} = X' \cap N_i$ exhibits a rational curve through $[X, x_1, x_2, x_3, x_4, \eta]$.
\end{proof}

\bibliography{bibliography}
\bibliographystyle{amsalpha}

\end{document}